\documentclass[reqno,12pt]{amsart}
\usepackage{amsfonts,amsmath,amssymb,fullpage,array,rotate,graphics}

\theoremstyle{plain} \newtheorem{Thm}{Theorem}[section]
\theoremstyle{plain} \newtheorem{Cor}[Thm]{Corollary}
\theoremstyle{plain} \newtheorem{Prop}[Thm]{Proposition}
\theoremstyle{plain} \newtheorem{Lemma}[Thm]{Lemma}
\theoremstyle{plain} 
\theoremstyle{definition} \newtheorem{Def}[Thm]{Definition}
\theoremstyle{remark} \newtheorem{Rem}[Thm]{Remark}
\theoremstyle{definition} \newtheorem{Ex}[Thm]{Example}

\def\pf{\noindent {\it Proof.}\hskip 8pt}
\def\qed{\hfill{\setlength{\fboxrule}{0.8pt}\setlength{\fboxsep}{1mm}
\fbox{\null}} \vskip 10pt}

\newcommand{\clearemptydoublepage}
             {\newpage{\pagestyle{empty}\cleardoublepage}}

\newcommand{\thmlist}%
          {\renewcommand{\theenumi}{\alph{enumi}}
                  \renewcommand{\labelenumi}{{\rm (\theenumi)}}
                            }
\newcommand{\factlist}%
          {\renewcommand{\theenumi}{\arabic{enumi}}
                  \renewcommand{\labelenumi}{{\rm (\theenumi)}}
                            }
\newcommand{\factlisti}%
          {
                  
                            }
\newcommand{\condlisti}%
          {
                  \renewcommand{\labelenumi}{{\rm(\theenumii)}}
                            }
\newcommand{\condlist}%
          {\renewcommand{\theenumi}{\roman{enumi}}
                  \renewcommand{\labelenumi}{{\rm (\theenumi)}}
                            }

\def\frak{\mathfrak}

\newcommand{\conv}{\mathop{\rm{conv}}}
\newcommand{\id}{\mathop{\rm{id}}}

\newcommand{\supp}{\mathop{\rm{supp}}}

\newcommand{\Ker}{\mathop{\rm{ker}}}

\newcommand{\abs}[1]{\left|\/#1\/\right|}

\newcommand{\inner}[2]{\langle#1,#2\rangle}

\newcommand{\C}{\ensuremath{\mathbb C}}

\newcommand{\D}{\ensuremath{\mathbb D}}

\newcommand{\R}{\ensuremath{\mathbb R}}
\newcommand{\Z}{\ensuremath{\mathbb Z}}
\newcommand{\N}{\ensuremath{\mathbb N}}

\renewcommand{\l}{\lambda}
\renewcommand{\a}{\alpha}
\renewcommand{\b}{\beta}
\newcommand{\e}{\varepsilon}

\newcommand{\fa}{\mathfrak{a}}

\newcommand{\frakacs}{\frak a_{\scriptscriptstyle{\C}}^*}

\newcommand{\rootstheta}{\langle\Theta\rangle}
\newcommand{\posrootstheta}{\langle\Theta\rangle^+}

\newcommand{\pedtheta}[1]{{#1}_{\scriptscriptstyle{\Theta}}}

\newcommand{\pedempty}[1]{{#1}_{\scriptscriptstyle{\emptyset}}}

\newcommand{\pedPi}[1]{{#1}_{\scriptscriptstyle{\Pi}}}

\newcommand{\cciwtheta}{C^\infty_c(\pedtheta{\frak a})^{\pedtheta{W}}}
\newcommand{\dcXsuppW}{\mathcal E'(X;W)}
\newcommand{\cciXsuppW}{C^\infty_c(X;W)}

\newcommand{\cciXWchi}{C^\infty_c(X;\chi)}
\newcommand{\cciXW}{C^\infty_c(X)^W}
\newcommand{\dcXWchi}{\mathcal E'(X;\chi)}
\newcommand{\dcXW}{\mathcal E'(X)^W}
\newcommand{\ccisupptheta}{C^\infty_c(\pedtheta{\frak a};\pedtheta{W})}
\newcommand{\dcsupptheta}{\mathcal E'(\pedtheta{\frak a};\pedtheta{W})}

\newcommand{\cciwW}{C_c^\infty(\frak a)^{\widetilde{W}}}
\newcommand{\cciaW}{C_c^\infty(\frak a)^{W}}
\newcommand{\dcsuppW}{\mathcal E'(\frak a;W)}
\newcommand{\cciaWchi}{C_c^\infty(\frak a;\chi)}
\newcommand{\cciasuppW}{C_c^\infty(\frak a;W)}

\begin{document}

\makeatletter
\title[Differential operators with hyperplane singularities]
{Support properties and Holmgren's uniqueness theorem\\
for differential operators with hyperplane singularities}
\author{Gestur \'{O}lafsson}
\address{Department of Mathematics, Louisiana State University, Baton Rouge,
LA 70803, U.S.A.}
\email{olafsson@math.lsu.edu}
\author{Angela Pasquale}
\address{D\'epartement et Laboratoire de Math\'ematiques,
Universit\'e de Metz, F-45075 Metz, France.}
\email{pasquale@math.univ-metz.fr}
\thanks{The first author was supported by NSF grants DMS-0070607,
DMS-0139783 and DMS-0402068. Part of this research
was conducted when the second author
visited LSU in February 2003. She gratefully
acknowledge financial support from NSF and
the Louisiana Board of
Regents grant \textit{Visiting Experts in
Mathematics}.  The final version of the
article was prepared while both authors were visiting the Lorentz
Center at Leiden University. They would like to thank E.~Opdam,
M.~de~Jeu, S.~Hille, E.~Koelink, W.~Kosters, M.~Pevzner and
F.~Bakker for their invitation.}
\date{}

\begin{abstract}
Let $W$ be a finite Coxeter group acting linearly on $\R^n$. In
this article we study support properties of $W$-invariant partial
differential operator $D$ on $\R^n$ with real analytic
coefficients. Our assumption is that the principal symbol of $D$
has a special form, related to the root system corresponding to
$W$. In particular the zeros of the principal symbol are supposed
to be located on hyperplanes fixed by reflections in $W$. We show
that $\conv (\supp Df) = \conv(\supp f)$ holds for all compactly
supported smooth functions $f$ so that $\conv(\supp f)$ is
$W$-invariant. The main tools in the proof are Holmgren's
uniqueness theorem and some elementary convex geometry. Several
examples and applications linked to the theory of special
functions associated with root systems are presented.
\end{abstract}
\maketitle
\makeatother

\section*{Introduction}
\noindent
Let $D$ be a linear partial differential operator on $\R^n$ with
constant coefficients. Then
a classical theorem of Lions and Titchmarsch states that, for every
distribution $u$ on $\R^n$
with compact support, the convex hulls of the supports of $Du$ and $u$
are equal:
\begin{equation} \label{eq:supportconstantPDO}
\conv(\supp Du)=\conv(\supp u), \qquad u \in \mathcal E'(\R^n).
\end{equation}
By regularization, this is equivalent to stating that for every
compactly supported smooth function $f$ on $\R^n$, the convex hulls of the
supports of $Df$ and $f$ are equal:
\begin{equation}
\conv(\supp Df)=\conv(\supp f), \qquad f \in C^\infty_c(\R^n).
\end{equation}
In fact, Lions \cite{Lio51} proved a more general version of the
support theorem, namely
\begin{equation} \label{eq:suppconv}
\conv\big(\supp(v * u)\big)=\conv(\supp v)+\conv(\supp u),
\qquad \forall\, v,\, u \in \mathcal E'(\R^n).
\end{equation}
We refer to \cite{Bo85} for an elementary proof of this theorem.
The first version (\ref{eq:supportconstantPDO}) obviously follows
from the third one (\ref{eq:suppconv}) by taking $u=\delta_0$, the
delta distribution at the origin.

The comparison of the supports of $f$ and $Df$ plays an important
role at several places in analysis. One typical situation is the
study of solvability of differential operators. Recall that a
linear partial differential operator
\begin{equation}\label{eq:D}
D=\sum_{\abs{I}\leq m} a_I(x) \partial^I
\end{equation}
with smooth coefficients $a_I$ is said to be {\em solvable in
$\R^n$}, provided $DC^\infty(\R^n)=C^\infty(\R^n)$, that is, if
for every $g \in C^\infty(\R^n)$ the differential equation $Df=g$
has a solution $f \in C^\infty(\R^n)$. The following theorem gives
a necessary and sufficient condition for the solvability of $D$
(see \cite{Treves}, Theorem 3.3).
\begin{Thm}\label{thm:solvable}
Let $D$ be a linear partial differential operator 
with smooth coefficients in $\R^n$. Then $D$
is solvable if and only if
the following two condition are satisfied:
\begin{enumerate}
\thmlist
\item {\em (semi-global solvability)} For every compact subset $K$ of
$\R^n$ and for every
$g \in C^\infty(\R^n)$ there is a function $f \in C^\infty(\R^n)$ so that
$Df=g$ on $K$.
\item {\em ($D$-convexity of $\R^n$)} For every compact subset $K$ of
$\R^n$ there is a compact
set $K'$ so that for every $f \in \mathcal C^\infty_c(\R^n)$ the inclusion
$\supp(D^t f) \subseteq K$ implies  $\supp f \subseteq K'$. 
Here $D^t$ denotes the
formal transpose of $D$.
\end{enumerate}
\end{Thm}

The support theorem of Lions and Titchmarsch implies that, for
every linear partial differential operator $D$ with constant 
coefficients, $\R^n$ is
$D$-convex, and that we can in fact take $K'=\conv K$. Observe
that, in this case, condition (b) also implies that the operator
$D^t$ is injective on $C^\infty(\R^n)$.

The main result of this paper is Theorem \ref{thm:maindistr},
which provides an extension of the support theorem of Lions and
Titchmarsch to specific (but very large) classes of invariant
singular linear partial differential operators with real analytic
coefficients and to distributions having a compact support with
invariant convex hull. The invariance considered here is with
respect to certain finite groups of orthogonal transformations
generated by reflections. In the 2-dimensional case, examples of
such groups are the groups of symmetries of regular $n$-agons. The
principal symbols of the examined differential operators are
allowed to vanish, but only in a precise way, along the reflecting
hyperplanes. See formula (\ref{eq:symbol}) below.

Some restrictions in generalizing the theorem of Titchmarsch and Lions
are of course needed. The following easy example shows that the theorem
cannot hold for arbitrary linear differential operators with
variable coefficients, even in the one-dimensional case and with real
analytic coefficients.

\begin{Ex} \label{ex:xdx}
Consider the differential operator $D=x \frac{d}{dx}$ on $\R$\,.
Let $u=\chi_{[0,1]}$ be the characteristic function of the interval
$[0,1]$\,. Then
$Du=-\delta_1$, where $\delta_1$ denotes the delta measure at $1$.
Therefore $\supp(Du)=\{1\}$ and $\supp u=[0,1]$ are convex and different.
\end{Ex}

In the one-dimensional case, the differential operator
$D=x \frac{d}{dx}$ of Example \ref{ex:xdx} belongs to the class of
differential operators to which our support theorem applies. Observe
that $D$ is an even differential operator which is singular at $x=0$.
For this very special differential operator, our theorem states that
$\conv(\supp Du)=\conv(\supp u)$ if $u \in \mathcal E'(\R)$
satisfies one of the following conditions:
\begin{enumerate}
\item $\supp u \subset ]0,+\infty[$\,,
\item $\supp u \subset ]-\infty,0[$\,,
\item $\conv(\supp u)$ is symmetric with respect to the origin $0$.
\end{enumerate}
Of course the distribution $\chi_{[0,1]}$ from Example \ref{ex:xdx}
 does not fulfill any of these conditions.

The core of the proof of Theorem \ref{thm:maindistr} is to show
that the considered situation allows us to apply Holmgren's
uniqueness theorem for comparing the size of the supports of $f$
and $Df$ when $f$ is a compactly supported smooth function with
the property that $\conv(\supp f)$ is invariant. The employ of
Holmgren's theorem is the reason for imposing to the coefficients
of the considered differential operators to be real analytic.

Several other authors employed Holmgren's uniqueness theorem for
proving $D$-convexity properties (see e.g. \cite{RW},
\cite{Chang}, \cite{vdBSconv}). Remarkable is nevertheless, that
the proof of our theorem is very elementary. It requires only
basic facts from convex geometry and an application of equation
(\ref{eq:suppconv}).

The article consists of two parts. The first part contains the
proof of Theorem \ref{thm:maindistr}, which does not require any
knowledge of symmetric spaces. In the second part, we discuss
several applications of our main results, including hypergeometric
differential operators, Bessel differential operators, shift
operators, Hamiltonian systems, and invariant differential
operators on symmetric spaces.

The solvability of $G$-invariant differential operators is one of
the fundamental problems in the analysis on a symmetric space
$G/H$ (see \cite{He2}, p. 275). Recall that Theorem
\ref{thm:solvable} holds more generally if $\R^n$ is replaced with
a 2nd countable smooth manifold (see \cite{Treves}, p. 14). Since
$G/H$ is a second countable smooth manifold, one obtains the
equivalence between global and semi-global solvability, provided
one can prove that $G/H$ is $D$-convex. The solvability of
invariant differential operators on Riemannian symmetric space was
proved by Helgason \cite{He73}. In the general pseudo-Riemannian
symmetric case, van den Ban and Schlichtkrull \cite{vdBSconv}
determined a sufficient condition for a $G$-invariant differential
operator $D$ ensuring that $G/H$ is $D$-convex. This condition
involves the degree of the polynomial which is the image of $D$
under the Harish-Chandra isomorphism. It is for instance always
satisfied when $G/H$ is split, i.e., has a vectorial Cartan
subspace. The Riemannian symmetric spaces of the noncompact type
are examples of split symmetric spaces and so are the $K_\epsilon$
space of Oshima and Sekiguchi, \cite{OS}.

As a first application of Theorem \ref{thm:maindistr}, we deduce in Section
\ref{section:appl} the $D$-convexity of Riemannian symmetric spaces $G/K$ of
noncompact type when $D$ is a $G$-invariant differential operator on $G/K$.
Our method, which is based on taking the radial component of $D$ along
the Cartan subgroup, is different from those of used in
\cite{He73} and \cite{vdBSconv}.

Support theorems play also an important role in harmonic analysis,
in particular in Paley-Wiener type theorems. These theorems
characterize the space of functions which are image, under a
suitable generalization of the Fourier transform, of the compactly
supported smooth functions. Applications in this setting appeared
first in the work of van den Ban and Schlichtkrull on Fourier
transforms on pseudo-Riemannian symmetric spaces \cite{BS}. The
basic idea, which we shall outline more precisely in Section
\ref{section:appl}, is to cancel the singularities appearing in a
wave packet $f$ by applying a suitable differential operator $D$.
The problem is to compare the size of the support of $Df$, which
can be easily determined, with the --hard to determine-- size of
the support of the original wave packet $f$. In fact, our need for
support properties like those stated in the present paper turned
up in the proof of a Paley-Wiener type theorem for the
$\Theta$-hypergeometric transform, which is a Fourier type
transform related to the theory of hypergeometric functions
associated with root systems. See \cite{OP4}. However, we point
out that Theorem \ref{thm:maindistr} is stated in a very general
setting and applies to many different situations. See Example
\ref{ex:diffrankone}, and the
Subsections \ref{ex:hyper}, \ref{ex:bessel} and
\ref{ex:shift} for several special cases.



Our paper is organized as follows. In Section
\ref{section:setting} we introduce the general setting in which
our extension of the theorem of Lions and Titchmarsch will be
proved. The main results, Theorems \ref{thm:maindistr} and
\ref{thm:mainfunction}, will be stated and proved in Section
\ref{section:thmsupport}. Section \ref{section:examples} presents
several concrete situations where our support theorem applies. The
presented examples are related to the theory of special functions
associated with root systems. The last section is devoted to
applications. We deduce the $D$-convexity of Riemannian symmetric
spaces of noncompact type when $D$ is an invariant differential
operator. Moreover, we describe how to employ Theorem
\ref{thm:mainfunction} for proving of Paley-Wiener type theorems
in the harmonic analysis on symmetric spaces and on root systems.

\section{Notation and setting} \label{section:setting}
\noindent
\subsection{Finite Coxeter Groups}
In this section we introduce the notation
and set up that will be used
in this article. In particular we introduce
the class of differential operators that will be
considered in this article and give few examples.

In the following $\fa$ stands for an real Euclidean vector space
of dimension $n$, i.e., $\fa\simeq \R^n$. Furthermore
$D$ will stand for a differential operator on $\fa$
with analytic coefficients.  We denote by
$\inner{\cdot}{\cdot}$ a (positive definite)
inner product on $\frak a$. Set $\abs{x}=\sqrt{\inner{x}{x}}$. For $\e>0$
we denote by $B_\e:=\{x \in \frak a: \abs{x}\leq \e\}$ the
closed Euclidean ball in $\frak a$ with center $0$ and radius $\e$.

Let $\frak a^*$ denote the real dual of
$\frak a$.
For each $\a \in \frak a^*\setminus\{0\}$ we denote by $y_\a$ the 
unique element of $\frak a$ satisfying
$\a(x)=\inner{x}{y_\a}$ for all $x \in \frak a$.
We set 
\begin{equation}
  \label{eq:xa}
  x_\a:=2y_\a/\inner{y_\a}{y_\a}
\end{equation}
and notice that $x_\a$ is independent of
$\inner{\cdot}{\cdot}$.
With each $\a \in \frak a^*\setminus \{0\}$ 
we associate the reflection $r_\a$ of 
$\frak a$ across the hyperplane $\mathcal H_\a:=\ker \a$. Thus
\begin{equation}
r_\a(x)=x-\a(x)x_\a, \qquad x\in \frak a.
\end{equation}
A finite set $\Delta \subset \frak a^* \setminus \{0\}$ is called a 
\emph{(reduced) root system} if the following conditions
holds for $\Delta$:
\begin{enumerate}
\renewcommand{\labelenumi}{(R\theenumi)}
\item If $\a \in \Delta$, then $\Delta \cap \R \alpha =\{\pm \alpha\}$;
\item If $\alpha,\beta \in \Delta$, then $r_\a(\beta) \in \Delta$.
\end{enumerate}
The elements of $\Delta$ are called \emph{roots}.
Observe that \textit{we are not requiring} that $\Delta$ contains a
basis of $\frak a^*$.  In particular, our definition allows $\Delta$ to be
the empty set.

A subset $\Pi$ of $\Delta$ is called a \emph{simple system} if
$\Pi$ is linearly independent and if any root in
$\Delta$ can be written as a linear combination of elements in $\Pi$
in which all non-zero coefficients are either all positive
or all negative. If $\Delta=\emptyset$, then
we set $\Pi=\emptyset$.

Let a simple system $\Pi$ of $\Delta$ be fixed. 
Set $\Delta^+:= \R_{\ge 0}\Pi \cap \Delta$,
where $\R_+=\{x\in \R\, :\, x\ge 0\}$. The elements
in $\Delta^+$ are said to be
positive.  Observe that
$-\a=r_\a(\a) \in \Delta$ for all $\a \in
\Delta$. Hence $\Delta=\Delta^+ \dot{\cup} (-\Delta^+)$.

Let $W \subset {\rm GL}(\frak a)$ be the 
group of orthogonal transformations of $\frak a$ generated by 
the reflections
$\{r_\a: \a \in \Delta\}$. It coincides with the group generated
by $\{r_\a\, :\, \a\in\Pi\}$. We set $W=\{\mathrm{id}\}$, if $\Delta$
is empty.  If $\dim \frak a=1$ and $\Delta\neq \emptyset$, then
 $W=\{\pm \id\}$.
The group $W$ is a \emph{finite Coxeter group}. Conversely, every finite 
Coxeter group originates from a root system as above. See e.g. \cite{GP}, 
Chapter 1. Among the finite Coxeter groups we find the \emph{Weyl groups} 
(for instance the finite groups of permutations) and the \emph{dihedral groups} 
(that is, the groups of symmetries of 
the regular $n$-gons).

The group $W$ acts on functions $f:\frak a \to \C$ according to 
\begin{equation} \label{eq:Wonf}
  (w \cdot f)(x):=f(w^{-1}x), \qquad w\in W,\, x \in \frak a\,.
\end{equation}
It also acts on compactly supported 
distributions $u \in \mathcal E'(\frak a)$ and on
differential operators $D$ on $\frak a$ by:
\begin{align}
\label{eq:Wonu}
<(w \cdot u),f>&:=<u,w^{-1} \cdot f>, \qquad w\in W,\,
f \in C^\infty(\frak a)\,,\\
\label{eq:WonD}
  (w \cdot D)f&:=w \cdot D(w^{-1} \cdot f), \qquad w\in W,\, 
f \in C^\infty(\frak a)\,,
  \end{align}
where we have written the pairing between distributions and functions
by $<u,f>:=u(f)$.
The function $f$ (resp. the compactly supported distribution $u$ or the 
differential operator $D$) is said to 
be \emph{$W$-invariant} provided $w\cdot f=f$ for all 
$w \in W$ 
(resp. $w \cdot u=u$ or $w\cdot D=D$ for all $w \in W$). 
For instance, if $\dim \frak a=1$ and $\Delta \neq \emptyset$, then 
$W$-invariant means even.
More generally, let $\chi$ be a character of $W$. 
Then $f$ is said \emph{to transform under $W$ according to $\chi$} if $w\cdot f=\chi(w) f$ for all
$w \in W$. This definition extends similarly to distributions
and differential operators.
Notice that if $f$ transforms under $W$ according to
a character $\chi$, then the
support of $f$ is $W$-invariant.

\subsection{$X$-elliptic polynomials and $X$-elliptic differential operators}
We will be studying differential operators with leading symbol of
a specific form. We will therefore need the following
definition:
\begin{Def}
Let $X$ be a non-empty $W$-invariant convex open subset of $\frak a$,
and
$P:\frak a\times \frak a^* \to \C$ a polynomial function.
We say that $P$ is a \textit{homogeneous $X$-elliptic polynomial}
if the following holds:
\begin{enumerate}
\renewcommand{\labelenumi}{(P\theenumi)}
\item
$P$ is a homogeneous polynomial in $\l \in \frak a^*$ with real analytic 
coefficients on $\frak a$, i.e. there is $m \in \N_0$ so that   
\begin{equation*}
P(x,\l)=\sum_{\abs{I}=m} a_I(x) \l^I\, 
\end{equation*}     
where $a_I(x)$ is real analytic,
and $\abs{I}:=\sum_{k=1}^n i_k$,
if $I=(i_1,\ldots ,i_n)\in\N_0^n$ is a multi-index.
\item
If $(x,\l) \in X \times \frak a^*$ and  $\l \neq 0$, then
$P(x,\l)\neq 0$.
\end{enumerate}
Let $D$ be a $W$-invariant linear partial differential operator 
on $\frak a$ with real analytic coefficients. We say that $D$ is 
\emph{$X$-elliptic} if its principal symbol is of the form
\begin{equation}\label{eq:symbol}
\sigma(D)(x,\l)=\left( p(\l) \; \prod_{\a \in \Delta} [\a(x)]^{n(\a)}
\right) P(x,\l) , \qquad
(x,\l)\in \frak a \times \frak a^*,
\end{equation}  
where $p(\l)$ is a homogeneous polynomial, 
$n(\a) \in \N_0:=\{0,1,2,\dots\}$ for all $\a \in \Delta$, and 
$P(x,\l)$
is a homogeneous $X$-elliptic polynomial. 
\end{Def}

In this paper we consider $X$-elliptic linear partial differential operators
as operators acting on functions or distributions on $X$.
 Note that the decomposition
of $\sigma (D)(x,\l )$ in (\ref{eq:symbol}) is in general not unique: if
$\alpha (x)\not= 0$ for all $x\in X$, then we can replace $P$ by
$\alpha (x)^k P(x,\lambda )$, $0<k\le n(\alpha)$, and replace
$n(\alpha)$ by $n(\alpha )-k$.

The class of $X$-ellipltic 
contains all elliptic partial differential operators on
$\fa$ by taking $X=\fa$ and $\Delta=\emptyset$.
But, more generally, the principal symbols of the considered operators
are allowed to vanish at the zeros of $p(\l)$ as well as along the hyperplanes
$\mathcal H_\a$ with $\a \in \Delta$.
Observe that the condition in (\ref{eq:symbol}) imposes a restriction only
on the principal part of the differential operators. In particular, suppose
$D_1$, $D_2$ are $W$-invariant linear partial differential operators with real analytic coefficients 
so that $\deg D_1 > \deg D_2$. If $D_1$ satisfies condition (\ref{eq:symbol}),
then the same is true for $D_1+D_2$.

\begin{Ex}
Any partial differential operator with constant coefficients 
$p(D)$ satisfies (\ref{eq:symbol}) when 
we choose $X=\frak a$ and $\Delta=\emptyset$. Indeed, in this case we do not 
impose any symmetry condition, and $\sigma(p(D))(\l)=p(\l)$ is of the form 
(\ref{eq:symbol}). 
\end{Ex}

\begin{Ex}[The one-dimensional case] 
\label{ex:diffrankone}
Suppose $\dim\frak a=1$ and $\Delta\neq \emptyset$. We shall identify 
$\frak a\equiv \frak a^*$ with $\R$. The possible subsets $X$ are the open 
intervals of the form $]-R,R[$ with $0< R \leq +\infty$. 
Then the differential operators considered are the even ordinary 
differential operators with real 
analytic coefficients and with principal symbol 
\begin{equation*}
\sigma(D)(x,\l)=x^n p(\l) P(x,\l).
\end{equation*}
Here, $n \in \N_0$, the polynomial $p(\l)$ is homogeneous,
and $P(x,\l)$ is homogeneous $X$-elliptic.
For $X=\R$, examples are $x\frac{d}{dx}$ and the Euler operator 
\begin{equation*}
E=x^2\frac{d^2}{dx^2}+a\, x\frac{d}{dx} +b\,, 
\end{equation*}
where $a,b \in \C$. Other examples, which are relevant 
in harmonic analysis, are constructed from the Jacobi and the Bessel 
differential operators on $]0,+\infty[$. Let $L$ and $L_0$ 
respectively denote the 
Jacobi and the Bessel differential operators, i.e.
\begin{align*}
L&=\frac{d^2}{dx^2}+ \Big[a \coth(x)+b \coth(2x)\Big] \frac{d}{dx}\,, \\
L_0&=\frac{d^2}{dx^2}+ \Big[a \frac{1}{x}+b \frac{1}{2x}\Big] 
\frac{d}{dx}\,. 
\end{align*} 
Then, for specific positive integral values of the
constants $a$ and $b$, the operator $L$ 
agrees with the radial part of the Laplace-Beltrami operator 
of a rank-one Riemannian symmetric spaces $G/K$ of noncompact type, 
whereas $L_0$ is the radial part of the Laplace-Beltrami operator 
of the corresponding rank-one Riemannian symmetric space $G_0/K$ 
of Euclidean type. 
Examples of such spaces $G/K$ are the (real, complex or quaternionic) 
hyperboloids; examples of spaces $G_0/K$ are given by $\R^n$
regarded as homogeneous space $M(n)/O(n)$. 
We refer to \cite{He2}, Ch. II, Propositions 3.9 and 3.13, for more 
information on radial parts of the Laplace-Beltrami operators.

The operators $D:=\sinh^2 x\cdot L$ and $D_0:=x^2\cdot L_0$ satisfy our
assumptions since they are even operators
on $\R$ with real analytic coefficients, and 
$\sigma(D)(x,\l)=x^2 \l^2 \big(\frac{\sinh x}{x}\big)^2$ and
$\sigma(D_0)(x,\l)=x^2 \l^2$.
Generalizations in more variables of these examples will be treated in
Section \ref{section:examples}. Note that, in the first case, we
can take $P(x,\lambda )=\big(\frac{\sinh x}{x}\big)^2\lambda^2$
or  $P(x,\lambda )=\big(\frac{\sinh x}{x}\big)^2$.
\end{Ex}

\begin{Ex}[The Calogero model] 
\label{ex:Calogero}
The Calogero model is a non-relativistic quantum mechanical system of 
$n+1$ identical particles on a line interacting pairwise. Such a system is 
described by the Hamiltonian 
\begin{equation*}
H_{\textrm{Cal}}(x)=-\frac{1}{2}\;\sum_{j=1}^{n+1} p^2_j + 
g^2 \sum_{1\leq i<j\leq n+1} \frac{1}{(x_i- x_j)^2}\,,\qquad (p_1,\ldots ,p_{n+1}),(x_1,\ldots ,x_{n+1}) \in
\R^{n+1},
\end{equation*}
where the positive constant $g^2$ is the coupling coefficient. 
The associated Schr\"odinger operator is
\begin{equation*}
S_{\textrm{Cal}}=-\frac{1}{2}\;\sum_{j=1}^{n+1} \partial^2_j +
g^2 \sum_{1\leq i<j\leq n+1} \frac{1}{(x_i- x_j)^2}\,.
\end{equation*}
See e.g. \cite{Pe90}, (3.1.1) and (3.1.14), I. See also \cite{Ca71}.
Let $\{e_1,\dots,e_{n+1}\}$ be the
standard basis of $\R^{n+1}$ and let $\inner{\cdot}{\cdot}$ be the standard
inner product on $\R^n$.
For $1\leq i, j\leq n+1$ with $i\neq j$ set $\a_{i,j}:=e_i-e_j$.
Then $\Delta:=\{\a_{i,j}:1 \leq i,j\leq n+1,
i \neq j\}$ is a root system of type $A_n$.
We take  $\Delta^+:=\{\a_{i,j}: 1\leq i<j\leq n+1\}$
as a set of positive roots. The
corresponding system of simple roots is $\Pi=\{\a_{j,j+1}:j=1,\dots,n\}$.
The finite Coxeter group $W$ associated to $A_n$ is the group
$\mathfrak{S}_n$ of permutations of the set $\{e_1,\dots,e_{n+1}\}$.
We can now write
\begin{equation*}
S_{\textrm{Cal}}(x)=-\frac{1}{2}\, L_{\frak a} + g^2 \sum_{\a \in \Delta^+}
\frac{1}{\inner{\a}{x}^2}
\end{equation*}
where $L_{\frak a}$ is the Euclidean Laplace operator on $\frak a$.
It follows that $S_{\textrm{Cal}}$ is $W$-invariant. If
$\pi(x):=\prod_{\a\in\Delta} \inner{\a}{x}$, then the  differential operator
$\pi(x) S_{\textrm{Cal}}(x)$ satisfies our requirements in (\ref{eq:symbol}).
This example will be generalized in Section \ref{section:examples}.
\end{Ex}

\begin{Rem}
Multiplication by $\sinh^2 x$ and $x^2$ in Example \ref{ex:diffrankone}
(resp. by $\pi$ in Example \ref{ex:Calogero}), which is used
to obtain $W$-invariant differential operators with analytic coefficients
in $\frak a$, is in fact inessential when dealing with even (resp. $W$-invariant functions).
Indeed, despite the singularity of their coefficients, the operators $L$ and $L_0$ map even
smooth functions on $\R$ in smooth functions. Likewise, the operator $S_{\textrm{Cal}}$ maps
smooth $W$-invariant functions on $\frak a$ into smooth functions. This remark will be stated,
in a more general form, in Section \ref{section:examples}.
\end{Rem}

\subsection{Function spaces}
We now introduce the class of functions that will be considered  
in this paper. As above, let $X$ be a fixed $W$-invariant open
convex subset of $\frak a$. Let $\dcXsuppW$ denote the space consisting 
of distributions $u$ on $\frak a$ so that $\conv(\supp u)$ is a 
$W$-invariant compact subset of $X$, and let $\cciXsuppW$ be the subspace 
of $\dcXsuppW$ consisting of $C^\infty$ functions.
Important subspaces of $\cciXsuppW$ are the spaces
$\cciXWchi$, formed by the smooth compactly supported functions 
$f: X \to \C$ which transform under $W$
according to the character $\chi$ of $W$. 
For instance, if $\chi$ is the trivial character, then we obtain
the space $\cciXW$ of $W$-invariant smooth functions on 
$X$ with compact support.
Similarly, inside $\dcXsuppW$ we find the spaces $\dcXWchi$ of 
compactly supported distributions on $X$ which 
transform according to $\chi$, and the space
$\dcXW$ of
$W$-invariant distributions on $W$ with 
compact support.

\begin{Ex}
If $\dim \frak a=1$ and $X=\frak a$, then we can identify
$X=\frak a\equiv \frak a^*$ 
with $\R$. Suppose first $\Delta=\emptyset$.
Then $W$ is trivial, and $\cciasuppW$ (resp. $\dcsuppW$) reduces to the 
space of $C^\infty$ functions (resp. distributions) on $\R$ with compact
support. If $\Delta\neq \emptyset$, then $W=\{\pm \id\}$.
In this case, $\cciasuppW$ (resp. $\dcsuppW$) is the space of smooth functions 
(resp. distributions) $f$ on $\R$ so that $\conv(\supp f)$ is a 
bounded interval of the form $[-R,R]$ for some $R>0$. 
The subspace $\cciaW$ consists of the compactly supported 
even smooth functions on $\R$; with $\chi$ equal to 
the sign character, the subspace $\cciaWchi$
consists of the compactly supported odd smooth functions on $\R$.
\end{Ex}

\section{The support theorem}
\label{section:thmsupport}
\noindent
In this section we proof the main result
of this article. This is the following version of the
the support theorem of Lions and Titchmarsch:

\begin{Thm}[The $W$-Invariant Support Theorem] \label{thm:maindistr}
Let $\emptyset\not= X \subseteq \frak a$ be open, convex and $W$-invariant.
Suppose $D$ is a $W$-invariant linear partial differential operator
 on $\frak a$ with real 
analytic coefficients and with principal symbol of the form
(\ref{eq:symbol}). Then 
\begin{equation*}
\conv\big(\supp Du\big)=\conv(\supp u)
\end{equation*} 
 for each $u \in \dcXsuppW$. 
\end{Thm}

We shall prove below that, by regularization, Theorem \ref{thm:maindistr}
is  equivalent to the following smooth version:

\begin{Thm} \label{thm:mainfunction}
Let $X$ and $D$ be as in Theorem \ref{thm:maindistr}.  Then 
\begin{equation*}
\conv\big(\supp Df\big)=\conv(\supp f)
\end{equation*} 
 for each $f \in \cciXsuppW$. 

Equivalently, for each $f \in \cciXsuppW$ and for every compact convex 
$W$-invariant
subset $C \subset X$ we have 
\begin{equation*}
\supp Df\subseteq C \quad \iff \qquad  \supp f \subseteq C\,.
\end{equation*} 
\end{Thm}

\smallskip
\noindent\textit{Proof of the equivalence of Theorems \ref{thm:maindistr} and
\ref{thm:mainfunction}.\;} 
It is clear that  Theorem \ref{thm:maindistr} implies
Theorem \ref{thm:mainfunction}.
For the other direction,
let $u \in \dcXsuppW$. Let $\{\psi_\e: \e>0\}$ be an
approximate identity with $\supp \psi_\e \subseteq B_\e$ for all $\e>0$. 
Since $X$ is a $W$-invariant, we
can, if necessary, replace $\psi_\e$ by
$\frac{1}{\abs{W}} \sum_{w \in W} (w \cdot \psi_\e)$.
We can therefore suppose
that $\psi_\e$ is $W$-invariant. It follows, in particular, that
$\supp \psi_\e$ and hence $\conv(\supp \psi_\e)$ are $W$-invariant.
Thus, by (\ref{eq:suppconv}),
$$\conv(\supp (u \ast \psi_\e))=\conv(\supp u)+\conv(\supp \psi_\e)$$
is a
$W$-invariant compact subset of $\frak a$. Since $X$ is open, we 
can choose $\e_0>0$ so that $u\ast \psi_\e \in \cciXsuppW$ for
$0<\e<\e_0$. The functions $u\ast \psi_\e$ converge to
$u$ in $\mathcal E'(\frak a)$ as $\e\to 0$, and $D(u \ast \psi_\e)=
Du \ast \psi_\e$. It is easy to check (without invoking deeper 
support theorems) that there exists $\delta(\e)>0$, converging to 
$0$ as $\e\to 0$, so that
\begin{equation} \label{eq:suppconvcontains}
\supp u \subseteq \supp(u\ast \psi_\e) + B_{\delta(\e)}\,.
\end{equation}
Replacing $u$ with $Du$, we obtain
\begin{align*}
\supp(Du\ast \psi_\e) &\subseteq \supp Du + B_\e\,,\\
\supp Du &\subseteq \supp(Du\ast \psi_\e) + B_{\delta(\e)}\,.
\end{align*}
The same inclusions are preserved for the convex hulls. By Theorem
\ref{thm:mainfunction}, 
we have
$$\supp(u\ast \psi_\e)=\supp(D(u\ast \psi_\e))\, .$$
Hence,
\begin{align*}
  \conv(\supp u) &\subseteq \conv (\supp(u\ast \psi_\e)) + B_{\delta(\e)}\\
         &= \conv \big(\supp(Du\ast \psi_\e)\big) + B_{\delta(\e)}\\
         &\subseteq \conv (\supp Du) + B_{\e+\delta(\e)},
\end{align*}
which implies that $\conv(\supp u) \subseteq \conv (\supp Du)$. Similarly,
$\conv (\supp Du) \subseteq \conv (\supp u)$. This proves Theorem 
\ref{thm:maindistr}.
\qed

Before going into the details of the proof of Theorem
\ref{thm:mainfunction}, let us briefly explain the main ideas involved.
Since $\supp Df \subseteq \supp f$, it suffices to show that
if $C$ is a compact convex $W$-invariant
subset of $X$, then $\supp Df\subseteq C$ implies 
that $\supp f \subseteq C$.
This will be proved by contradiction. The main tool will be
Holmgren's uniqueness theorem.

\begin{Thm}[Holmgren's uniqueness theorem] 
  \label{thm:holmgren}
Let $\emptyset\not=\Omega \subseteq \R^n$ be open,
and let $\varphi$ be a real valued function in 
$C^1(\Omega)$. Let $D$  
be a linear partial differential operator with analytic coefficients 
defined in $\Omega$. Let $\sigma(D)$ 
denote the principal symbol of $D$.

Suppose that $x_0$ is a point in $\Omega$ such that
\begin{equation}
  \label{eq:nonvanishingsymb}
  \sigma(D)(x_0,d\varphi(x_0))\neq 0.
\end{equation}
Then there exists a neighborhood $\Omega' \subseteq \Omega$ of $x_0$ 
with the following property: If the
distribution $u\in \mathcal D'(\Omega)$
is annihilated by $D$, i.e., $Du=0$,
and $u$ vanishes on the set $\{x\in \Omega
: \varphi(x)>\varphi(x_0)\}$,
then $u$ must also vanish in $\Omega'$. 
\end{Thm}
\pf This is Theorem 5.3.1 in \cite{Hoermold}.
\qed
For non-elliptic differential operators the delicate matter is to 
choose points $x_0$ for which condition (\ref{eq:nonvanishingsymb}) is
fulfilled.  
This is exactly the kind of difficulty one encounters in the proof
of support theorems.
Coming back to the streamline ideas of the proof of 
Theorem \ref{thm:mainfunction}, let $C$ be as above. Set
$S:=\supp f$. Then $S \subseteq X$. To reach a contradiction, we 
assume that $S\not\subseteq C$. Let $x_0 \in S \setminus C$.
As $C$ is convex and compact, there exists a hyperplane strictly
separating $x_0$ and $C$, i.e., we
can find $\l_0\in\in \frak a^*$ such that
\begin{equation*}
 \max_{y \in C} \l_0(y) < \l_0(x_0). 
\end{equation*}
Without loss of generality we can also assume that
\begin{equation*}
  \l_0(x_0) = \max_{x \in S} \l_0(x)\, .
\end{equation*}
Otherwise  we translate the hyperplane to the boundary of $S$ in the
direction opposite to $C$. In this way, the entire set $S$ lies inside
a closed half-space supported by the hyperplane 
$$\mathcal H_0:=\{x \in \frak a: \l_0(x)=\l_0(x_0)\}\,.$$
Our plan is to apply Holmgren's uniqueness theorem to $\Omega=X
\setminus C$, $\varphi=\l_0$ and $u=f$. Note that in this case
 $d\varphi=\l_0$ is constant and non-zero.
Observe also that $Df\equiv 0$ on $\Omega$ and $f\equiv 0$
on the side of $\mathcal H_0$ not containing $C$ (which is described by the 
equation $\varphi(x)>\varphi(x_0)$). If the principal symbol of $D$ were
not zero at $(x_0,\l_0)$, then all assumptions would be satisfied, and we
could conclude that $f\equiv 0$ in a neighborhood of $x_0$. This would yield
the required contradiction because $x_0 \in S=\supp f$. 

Since $\l_0 \neq 0$, the condition
$\sigma(D)(x_0,\l_0)\neq 0$ is equivalent to
\begin{equation}\label{eq:nonzerosym}
  p(\l_0) \prod_{\a \in \Delta} [\a(x_0)]^{n(\a)}\neq 0.
\end{equation}
This might not be satisfied by the chosen pair $(x_0,\l_0)$.
It is even possible to have a situation where
there is no choice of $\l_0$ for which the above procedure could guarantee 
that $\a(x_0) \neq 0$. 
Figure 1 sketches an example in which this problem arises because of 
the ``corner'' at the boundary of $S$. Note that the set $S$ in this
example is also convex and invariant with respect to the group 
$W=\{\id, r_v\}$, where $r_v$ denotes the reflection with respect to the
$v$-axis. 

\begin{figure}
\parbox{7truecm}{\includegraphics{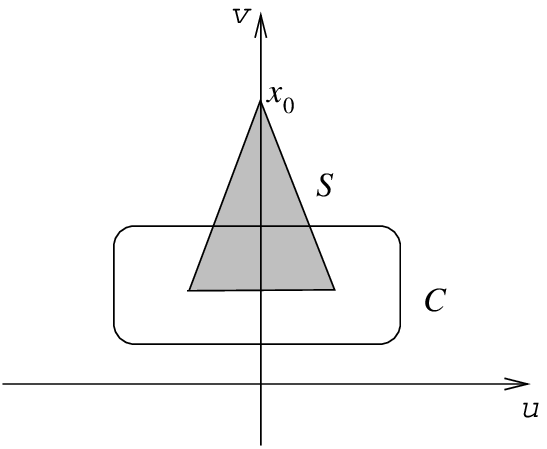}}
\parbox{5truecm}{$\frak a=X=\R^2$,\\
$\Delta=\{\pm (1,0)\}$ of type $A_1$, \\
$\Sigma^+=\Pi=\{(1,0)\}$}
\caption{}
\end{figure}

In case $ p(\l_0) \prod_{\a \in \Delta} [\a(x_0)]^{n(\a)}=0$,
the above procedure must be modified.
The first step is to show that it
suffices to consider the case of smooth functions 
$f$ with the property that the convex hull $C'$ of $S:=\supp f$ 
is $W$-invariant and has $C^1$ boundary. The point is that
the $x_0$, selected as above, will
always belong to the boundary
$\partial(C')$ of $C'$, and that $\mathcal H_0$ is a
supporting hyperplane for $C'$. For the modified procedure  
we need some preparations.

For $\e>0$ let $B_\e\subset \frak a$ denote the closed ball centered at
the origin and with radius $\e$. Let $B_\e(C):=C+B_\e$ be the
(closed) $\e$-neighborhood of $C$.

\begin{Lemma} \label{lemma:smoothbdry}
Let $C$ be a proper convex subset of $\frak a\equiv \R^n$ 
with nonempty interior, and let $\e>0$. 
Then $B_\e(C)$ is an $n$-dimensional convex subset of 
$\frak a$ and its boundary $\partial B_\e(C)$ is a $C^1$ 
$(n-1)$-dimensional submanifold of $\frak a$.
\end{Lemma}
\pf
This is Satz 17.6 in \cite{NGL}. \qed

\begin{Lemma} \label{lemma:withsmoothbdry}
Let $D$ be a differential operator on $\frak a$.
Suppose there exists a
smooth compactly supported function $\widetilde f$ on
$\fa$
and a compact convex subset $\widetilde C \subset \frak a$ with nonempty
interior such that
$$\supp D\widetilde f \subseteq \widetilde C \qquad\textrm{but}\qquad
\supp \widetilde f 
\not\subset \widetilde C.$$
Then there exists a function $f \in C_c^\infty(\frak a)$ and a compact convex
subset $C  \subset \frak a$ with nonempty interior such that
$$\supp Df \subset C \qquad\textrm{and} \qquad \supp f \not\subset C\, .$$
Moreover,
the boundary $\partial \big(\conv(\supp f)\big)$ of the convex hull of 
$\supp f$ is a $C^1$-manifold.

If $\widetilde{f}\in \cciXsuppW$, where $X$ is a $W$-invariant open convex  subset
of $\frak a$,
and  $\widetilde C \subseteq X$ is $W$-invariant.
Then we can choose $f \in \cciXsuppW$
and $C \subseteq X$ to be 
$W$-invariant.
\end{Lemma}
\pf 
Let $\{\psi_\e:\e>0\}$ be an approximate identity with 
$\supp \psi_\e = B_\e$ for all $\e$. 
Then 
$\widetilde f \ast \psi_\e \in C^\infty$ and  
\begin{equation*}
\supp\big(D(f\ast \psi_\e)\big)=\supp\big(Df\ast \psi_\e\big) 
\subseteq \supp Df + B_\e \subseteq B_\e(\widetilde C).
\end{equation*}
Notice that (\ref{eq:suppconv}) applied to compactly supported smooth 
functions implies that for all $g \in C_c^\infty$ we have
\begin{equation}
  \label{eq:supp}
 \supp g \subseteq B_\e\big(\conv(\supp g)\big) 
         =\conv(\supp g)+B_\e 
         =\conv\big(\supp(g\ast \psi_\e)\big).
\end{equation}
Hence, there exists $\e_0>0$ so that $\supp(\widetilde f \ast \psi_{\e_0})
\not\subseteq B_{\e_0}(\widetilde C)$. Otherwise  (\ref{eq:supp})
would imply that
$\supp \widetilde f \subseteq \conv\big(\supp(\widetilde
f\ast \psi_\e)\big)
    \subseteq B_{\e}(\widetilde C)$
for all $\e>0$, and hence $\supp \widetilde f \subseteq \widetilde C$.
As in (\ref{eq:supp}), we have
\begin{equation*}
  \conv\big(\supp(\widetilde f\ast \psi_{\e_0})\big)=
  B_{\e_0}\big(\conv(\supp \widetilde f)\big).
\end{equation*}
Therefore, by Lemma \ref{lemma:smoothbdry}, 
$\conv\big(\supp(\widetilde f\ast \psi_{\e_0})\big)$ has $C^1$-boundary. 
We can thus select
$f=\widetilde f\ast \psi_{\e_0}$ and $C=B_{\e_0}(\widetilde C)$.

Finally, suppose that $\conv(\supp \widetilde{f})$ and 
$\widetilde{C}$
are $W$-invariants subset of $X$. Since $X$ is open, then, 
by choosing a sufficiently small $\e_0>0$, 
so are also $\conv\big(\supp(\widetilde f\ast \psi_{\e_0})\big)=
  B_{\e_0}\big(\conv(\supp \widetilde f)\big)$ and  
 $B_{\e_0}(\widetilde{C})$.
 \qed

In the following we shall suppose that $f$ and $C$ are chosen as in Lemma
\ref{lemma:withsmoothbdry}. As before, we set $S:=\supp f$ and 
$C'=\conv(\supp f)$. We also fix $y_0 \in S \setminus C$.
We now proceed to the selection of the pair $(x_0,\l_0)$. 

\begin{Lemma}
Let $y_0 \in S\setminus C$ be fixed, and let
$$U_0:=\{\l\in \frak a^*: \max_{y \in C} \l(y)<\l(y_0)\}.$$
Then $U_0$ is a nonempty open subset of $\frak a^*$.
\end{Lemma}
\pf
The separation properties of compact convex sets ensure that $U_0 \neq
\emptyset$.
Observe that
\begin{align*}
U_0&=\{\l \in \frak a^*:\max_{y\in C} \l(y-y_0)<0\}\\
   &=\{\l \in \frak a^*:\min_{y\in y_0-C} \l(y)>0\}.
\end{align*}
Hence $U_0=h^{-1}(]0,+\infty[)$ where $h(\l):=\min_{y\in y_0 -C} \l(y)$.
This proves the lemma because $h$ is proper and convex, hence continuous (see
\cite{Rock}, Theorems 5.3 and 10.1).          
  \qed

Since $U_0$ is open and nonempty we can choose $\l_0\in U_0$ with the 
following properties:
\begin{enumerate}
\renewcommand{\labelenumi}{(\theenumi)}
\item $p(\l_0) \neq 0$,
\item $\inner{\l_0}{\a} \neq 0$ for all $\a \in \Delta$
(i.e. $\l_0(x_\a)\neq 0$ for all $\a \in \Delta$).
\end{enumerate}
Since $C'=\conv S$,
\begin{Lemma}[Choice of $(x_0,\lambda_0)$]\label{le:cofx0}
Let $\l_0$ be as above.
Then there exists $x_0 \in S$
such that $\l_0(x_0)=\max_{x\in C'} \l_0(x)$.
Furthermore $x_0 \in \partial S \cap \partial C'$ and $x_0 \in X$.
Finally
\begin{equation*}
  \l_0(x_0)=\max_{x\in S} \l_0(x) \geq \l_0(y_0) >  \max_{y \in C} \l_0(y).
\end{equation*}
\end{Lemma}
\begin{proof} This follow as $C'$ is the convex hull of $S$.\end{proof}

When $\Delta\neq \emptyset$, the $W$-invariance of the situation plays 
a role because of the following lemma.

\begin{Lemma} \label{lemma:tangent}
Suppose $C'$ is a $W$-invariant compact convex subset of
$\frak a$ with nonempty interior and $C^1$-boundary. Let
$\a\in\Delta$ and $x\in \partial C'$ such that $\a(x)=0$. Let
$T_x(\partial C')$ denote the tangent space to $\partial C'$ at $x$
(regarded as subspace of $\frak a$) and let $x_\a$ be as in (\ref{eq:xa}).
Then $x_\a \in T_x(\partial C')$.
\end{Lemma}

\pf
The reflection $r_\a$ of $\frak a$ across the hyperplane $\Ker \a$ 
maps $C'$ onto itself and
fixes $x$. Hence it maps $(C')^0$ (resp. $\partial C'$) onto itself. The
derived
involution $r_\a$ is therefore an automorphism of $T_x(\partial C')$. Hence
$r_\a(N)=\pm N$, where $N$ is the outer normal vector to $\partial C'$ 
at $x$.
Since $r_\a N\neq -N$ by invariance of $C'$ under $r_\a$, we conclude that
$r_\a N=N$, i.e. $N\in \Ker \a=x_a^\perp$. Thus $x_\a\in x_\a^{\perp\perp} 
\subseteq N^\perp =T_x(\partial C')$.        
\qed

\begin{Prop} \label{prop:nonzerosym}
Let $x_0$ and $\l_0$ be
as in Lemma \ref{le:cofx0}. Then  $(x_0,\l_0)$ satisfies (\ref{eq:nonzerosym}).
\end{Prop}
\pf
The element $\l_0$ has been chosen so that $p(\l_0) \neq 0$.  
This suffices to prove (\ref{eq:symbol}) when $\Delta =\emptyset$. 
If $\Delta\neq \emptyset$, then it remains to show that $\a(x_0)\neq 0$ 
for all $\a \in \Delta$.
Since $\l_0(x_0)=\max_{x\in C'} \l_0(x)$, the hyperplane 
$\mathcal H_0:=\{x \in\frak a: \l_0(x)=\l_0(x_0)\}$ is the
supporting hyperplane for $C':=\conv(\supp f)$ at $x_0$. The set
$C'$ is $W$-invariant and its boundary 
$\partial C'$ is $C^1$. Hence $\mathcal H=x_0+T_{x_0}\partial C'$. 
If $\a(x_0)=0$ for some $\a \in \Delta$, then $x_\a\in T_{x_0} 
\partial C'$ by Lemma \ref{lemma:tangent}, i.e. $x_0+x_\a\in \mathcal H_0$. 
Thus $\l_0(x_0)+\l_0(x_\a)=\l_0(x_0+x_\a)=\l_0(x_0)$, which implies 
$\l_0(x_\a)=0$, against our choice of $\l_0$. 
Thus $\a(x_0)\neq 0$ for all $\a \in
\Delta$.
\qed           

\smallskip
\noindent \textit{Proof of Theorem \ref{thm:mainfunction}.\;}
Arguing by contradiction, we assume that there exists 
$f\in \cciXsuppW$ and a $W$-invariant subset $C$ of $X$ 
as in Lemma \ref{lemma:withsmoothbdry}. We select $(x_0,\l_0)$ 
as in Proposition \ref{prop:nonzerosym}. Set $\Omega=X
\setminus C$ and $\varphi=\l_0$. Hence $\mathcal H_0:=\{x \in
\Omega: \l_0(x)=\l_0(x_0)\}$ is a supporting hyperplane for $C':=
\conv(\supp f)$, 
and $d\varphi=\l_0$ is constant and non-zero.
Moreover, $Df=0$ on $\Omega$ and $f\equiv 0$
on the side of $\mathcal H_0$ not containing $C$ (which is described by the 
equation $\l_0(x)>\l_0(x_0)$). Proposition \ref{prop:nonzerosym} ensures
that $\sigma(D)(x_0,\l_0)\neq 0$. Holmgren's Uniqueness Theorem then 
implies that $f\equiv 0$ in a neighborhood of $x_0$. This gives the
the required contradiction because $x_0 \in \supp f$. 
\qed

Before concluding this section we prove some immediate consequences of 
Theorem \ref{thm:maindistr}. Recall that the transpose of the
partial differential operator $D$ given by (\ref{eq:D}) is
\begin{equation}
  \label{eq:eq:Dt}
  D=\sum_{\abs{I}\leq m} (-1)^{\abs{I}} a_I(x) \partial^I\,.
\end{equation}
Hence $D^t$ belongs to the class of invariant differential operators
considered in this papers if so does $D$. In fact, the principal
symbols of $D$ and $D^t$ are linked by the relation
\begin{equation*}
  \sigma(D^t)(x,\l)=(-1)^m \sigma(D)(x,\l)\,, \qquad (x,\l) \in 
\frak a \times \frak a^*.
\end{equation*}

\begin{Cor}
  \label{cor:maindistr}
Let $D$ be a $W$-invariant differential operators $D$ on $\frak a$ with real 
analytic coefficients and with principal symbol of the form
(\ref{eq:symbol}). Then the following properties are true.
  \begin{enumerate}
\thmlist
  \item $D$ is injective on $\cciXsuppW$.
  \item For all $f \in \cciXsuppW$ we have
\begin{equation*}
\conv\big(\supp D^t f\big)=\conv(\supp f)\,.
\end{equation*} 
  \end{enumerate}
\end{Cor}

\section{Applications to some special differential operators} 
\label{section:examples}

\noindent
In this section we present some examples linked to the theory of
special functions associated with root systems. In these examples
the group $W$ is a parabolic subgroup of a fixed Coxeter group
$\widetilde{W}$ acting on $\frak a\simeq \R^n$, and
the differential operator $D$ is in fact
invariant under the larger finite Coxeter group $\widetilde{W}$.
The general situation corresponds to a (not necessarily reduced) root systems
$\Sigma$. A \emph{root system} is a finite set $\Sigma \subset \frak a^*
\setminus \{0\}$ satisfying condition (R2) of Section
\ref{section:setting}. In this section we will
also assume that
$\Sigma$ satisfies the following additional conditions:
\begin{enumerate}
\renewcommand{\labelenumi}{(R0)}
\item $\Sigma$ spans $\frak a^*$;
\renewcommand{\labelenumi}{(R2)}
\item $\Sigma$ is \emph{crystallographic}, that is 
\begin{equation*}
2\frac{\inner{\a}{\beta}}{\inner{\a}{\a}}\in \Z
\end{equation*}
for all $\a, \b \in \Sigma\,$.
\end{enumerate}

Crystallographic root systems arise naturally
in several places in algebra and analysis.
In particular, they are relevant in the
theory of real Lie algebras, Lie groups, and symmetric spaces.

If $\Sigma$ is a root system, then $\Delta:=\{\a \in \Sigma: 2\a
\notin \Sigma\}$ is a reduced root system according to the
definition of Section \ref{section:setting}. The finite Coxeter
group associated to $\Sigma$ is the
finite Coxeter group $\widetilde{W}$
associated to $\Delta$. It is also
called the \emph{Weyl group} of $\Sigma$.
A \emph{multiplicity function} is a $\widetilde{W}$-invariant function
$m:\Sigma \to \C$.
For $\a \in \Sigma$, we adopt the common notation $m_\a$ to denote $m(\a)$.

The set
$\Sigma^+$ of positive roots and the set $\Pi$ of positive simple
roots can be defined as in the case of reduced root systems.
Because of (R0), $\Pi$ is a basis of $\frak a^*$.
Fix a $\Pi$ set of simple roots in $\Sigma^+$.
For each subset $\Theta$ of $\Pi$ we define $\rootstheta$
as the set of elements of $\Sigma$ which are linear combinations of elements
from $\Theta$, i.e., $\rootstheta :=\Z \Theta \cap \Sigma$.
It is itself a root system, but in general it does not satisfy  (R0).
We denote the corresponding finite Coxeter group by $\pedtheta{W}$.
Note that $\pedtheta{W}\subseteq \widetilde{W}$ is
generated by
the reflections $r_\a$ with $\a \in \Theta$.
For instance, $\pedPi{W}=\widetilde{W}$
and $\pedempty{W}={\id}$. We also set $\posrootstheta:=\rootstheta \cap
\Sigma^+$ for the set of positive roots inside $\rootstheta$.
The Coxeter group $\pedtheta{W}$ will play the role of the group $W$ of the
previous sections. Recall that a subgroup $W$ of $\widetilde{W}$ is called
a \emph{parabolic subgroup} if it is of the form $\pedtheta{W}$ for
some $\Theta \subseteq \Pi$. The parabolic subgroups can also be characterized as
those subgroups of $\widetilde{W}$ that stabilize a subspace of $\frak a$. Thus
$W\subseteq \widetilde{W}$ is a parabolic subgroup if and only if there exists
a subspace $\frak b\subseteq \fa$ such that
$$W=\{w\in \widetilde{W} \, :\, w(\frak b )=\frak b\}\, .$$

A \emph{chamber} in $\frak a$ a connected component of
$\frak a \setminus \cup_{\a \in \Sigma} \mathcal H_\a$
(see \cite{Bou}, Ch. V, \S 3).
>From now on we fix the chamber
$\frak a^+:=\{x \in \frak a: \text{$\a(x)>0$ for all $\a \in \Pi$}\}$.
It is an open simplicial cone with vertex
$0$ (\emph{loc. cit.}, Ch. V, \S 3.9, Proposition 7 (iii)), and its closure
$\overline{\frak a^+}:=\{x \in \frak a: \text{$\a(x)\geq 0$
for all $\a \in \Pi$}\}$
is a fundamental domain for the action of $\widetilde{W}$ on $\frak a$
(\emph{loc. cit.}, Ch. V, \S 3.3, Theorem 2).

We define 
\begin{equation}\label{eq:aTheta}
  \pedtheta{\frak a}:=\pedtheta{W}(\overline{\frak a^+})^0,
\end{equation}
where $\null^0$ denotes the interior. 
For instance, $\pedempty{\frak a}=\frak a^+$ and $\pedPi{\frak a}=\frak a$.
In the following examples, the sets
$\pedtheta{\frak a}$ will play the role of the set $X$ appearing in Theorem
\ref{thm:maindistr}.

The set $\overline{\pedtheta{\frak a}}=\pedtheta{W}(\overline{\frak a^+})$
is the smallest cone 
in $\frak a$ which is closed, $\pedtheta{W}$-invariant and contains
$\overline{\frak a^+}$.
It is a union of closed chambers. 
The open polyhedral cone in $\frak a^*$
\begin{equation*}
  \pedtheta{C}:=\sum_{\a \in \Sigma^+ \setminus \posrootstheta} \R^+x_\a
\end{equation*}
has the dual cone
\begin{align*}
  \pedtheta{C}^*:&=\{x \in \frak a: \text{$\inner{x}{y} \geq 0$
for all $y \in \pedtheta{C}$}\}\\
             &=\{x \in \frak a: \text{$\a(x) \geq 0$  
           for all $\a \in \Sigma^+ \setminus \posrootstheta$}\}.  
\end{align*}
$\pedtheta{C}^*$ is a closed convex cone. It
is the  intersection of the closed hyperplanes defined by roots,
and hence a union of closed chambers in $\frak a$.

\begin{Lemma} \label{lemma:cone}
We have
\begin{equation*}
\overline{\pedtheta{\frak a}}=\pedtheta{C}^*.
\end{equation*}
Consequently, $\overline{\pedtheta{\frak a}}$ is a closed 
$\pedtheta{W}$-invariant 
convex cone in $\frak a$ and also its interior 
$\pedtheta{\frak a}$ is convex.  
Moreover, 
\begin{equation}\label{eq:pedthetaa}
\pedtheta{\frak a}=\{x \in \frak a: \text{$\a(x) > 0$  
           for all $\a \in \Sigma^+ \setminus \posrootstheta$}\}.
\end{equation}

\end{Lemma}
\pf
This was proven in \cite{Pa04}, Lemma 3.4, when $W$ is a Weyl group.
The same proof applies also to the more general case of finite Coxeter 
groups.
\qed

Specializing the notation of Section \ref{section:setting} to this context, 
we consider the space $\ccisupptheta$ of smooth functions 
$f: \frak a \to \C$ with the property that 
$\conv(\supp f)$ is a $\pedtheta{W}$-invariant compact subset
of $\pedtheta{\frak a}$. Its subspace of $\pedtheta{W}$-invariant functions
on $\pedtheta{\frak a}$ with compact support is $\cciwtheta$. 
Furthermore, $\dcsupptheta$ is the space  consisting 
of distributions $u$ on $\frak a$ so that $\conv(\supp u)$ is a 
$\pedtheta{W}$-invariant compact subset of
$\pedtheta{\frak a}$.


\subsection{Hypergeometric differential operators}\label{ex:hyper}
As before let $\frak a\simeq \R^n$ be a finite dimensional Euclidean space and let
$\Sigma$ be a (non-necessarily reduced) root system in $\frak a^*$.
Further, let $\widetilde{W}$ the be the corresponding Weyl group  and
let
$m$ be a multiplicity function on $\Sigma$.
Heckman and Opdam associated with such a triple $(\frak a, \Sigma, m)$
a commutative family of $\widetilde{W}$-invariant differential
operators on $\frak a$, the \emph{hypergeometric differential
operators}, having meromorphic coefficients of
a specific type. More precisely, their coefficients are meromorphic
functions on the complexification $\frak a_\C$ of $\frak a$ and their
singularities are cancelled by multiplication by a suitable power
of the Weyl denominator
\begin{equation} \label{eq:delta}
 \delta(x):=\prod_{\a\in\Sigma^+} \sinh \a(x)\,.
\end{equation}
We refer to \cite{HeckBou}, \cite{Opd00} and \cite{HS} for more information
on Heckman-Opdam's theory of hypergeometric differential operators.

For special values of $m$, the triple $(\frak a, \Sigma, m)$ arises from a
Riemannian symmetric space $G/K$ of the noncompact type, i.e., $G$
is a noncompact connected semisimple Lie group with
finite center and $K\subset G$ is a maximal
compact subgroup, see
Section \ref{section:appl}. In this case,
the hypergeometric differential operators
coincide with the radial parts (with respect
to the $K$-action) on $\frak a^+$
of the $G$-invariant differential operators on $G/K$.  Here we identify
$\frak a^+$ with its diffeomorphic image in $G/K$ under the exponential map
(usually denoted $A^+$ in the literature on analysis on symmetric spaces).
For instance, the hypergeometric differential operator
\begin{equation} \label{eq:L}
L:=L_{\frak a}+\sum_{\a \in \Sigma^+}m_\a \coth \a \;\partial(y_\a)
\end{equation}
coincides with the radial part of the Laplace-Beltrami operator of $G/K$.
Notice that $L$ is singular on the hyperplanes $\mathcal{H}_\alpha$, 
$\alpha\in \Sigma$.
Furthermore, $L$ is the multivariable analog of the Jacobi differential
operator of Example \ref{ex:diffrankone}.
In (\ref{eq:L}), $\partial(y)$ denotes the directional 
derivative in the direction of
$y \in \frak a$, and $L_{\frak a}$ is the Laplace operator on the 
Euclidean vector space
$\frak a$, that is, $L_{\frak a}=\sum_{j=1}^n \partial(x_j)^2$ where
$\{x_j\}_{j=1}^n$ is an orthonormal basis of $\frak a$.

Set
\begin{align*}
 \pedtheta{\pi}(x)&:=\prod_{\a\in\rootstheta^+} \a(x)\,,\\
 \pedtheta{\delta}(x)&:=\prod_{\a\in\rootstheta^+} \sinh \a(x)\,,\\
 \pedtheta{\delta}^c(x)&:=\prod_{\a\in\Sigma^+ \setminus
           \rootstheta^+} \sinh \a(x)\,,
\end{align*}
with the usual convention that empty products are equal to $1$. 
We write $\pi(x)$ instead of $\pedPi{\pi}(x)$.
Finally, define $D:=\delta^2 \cdot L$. The principal 
symbol of $D$ is
\begin{equation*}
  \sigma(D)(x,\l)=\delta(x)^2\cdot \sigma(L_{\frak a})(x,\l)
                 = \inner{\l}{\l}  \pedtheta{\pi}(x)^2
 P(\l,x)
\end{equation*}
where 
\begin{equation}\label{eq:Phyp}
  P(\l,x):=[\pedtheta{\delta}^c(x)]^2 
\Big[\frac{\pedtheta{\delta}(x)}{\pedtheta{\pi}(x)}\Big]^2.
\end{equation}
Because of Lemma \ref{lemma:cone}, each $\a \in \Sigma^+ \setminus
           \rootstheta^+$ is positive on $\pedtheta{\frak a}$.
Therefore $P(\l,x)$ is a homogeneous $\pedtheta{\frak a}$-elliptic
polynomial on $\frak a^* \times \frak a$ (of degree $0$ in $\l \in
\frak a^*$).
More generally, for each hypergeometric differential operator $D_0$ 
there is $k \in \N$ so that 
the linear partial 
differential operator $D:=\delta(x)^{2k} D_0$ is a $\widetilde{W}$-invariant
differential operator with real analytic coefficients and with
principal symbol of the form (\ref{eq:symbol}) where $P(\l,x)$ is
a homogeneous $\pedtheta{\frak a}$-elliptic polynomial.
In this case, we say that $D=\delta (x)^kD_0$ is
a \emph{regularization} of $D_0$.
The nature of the principal symbol of $D$ can in fact be deduced
by the explicit representation of the hypergeometric differential 
operators in terms of Cherednik operators (see e.g. \cite{HeckBou}).
According to this representation, each hypergeometric differential operator 
can be constructed as the ``differential part'' of the 
differential-reflection operator associated with a 
$\widetilde{W}$-invariant polynomial in the symmetric algebra 
over $\frak a_\C$ by means of Cherednik operators. If $D_0$ is 
associated with the $\widetilde{W}$-invariant polynomial $p$, 
then the principal part of the regularized 
operator $D=\delta(x)^{2k} \cdot D_0$ is
\begin{align*}
  \sigma(D)(x,\l)&= p_h(\l) \pi(x)^{2k}
    \Big[ \frac{\delta(x)}{\pi(x)}\Big]^{2k}\,\\
                  &= p_h(\l) \pedtheta{\pi}(x)^{2k} P(\l,x)
\end{align*}
where $p_h(\l)$ is the highest homogeneous part of $p(\l)$ and 
\begin{equation*}
 P(\l,x):=[\pedtheta{\delta}^c(x)]^{2k} 
\Big[\frac{\pedtheta{\delta}(x)}{\pedtheta{\pi}(x)}\Big]^{2k}.
\end{equation*}

In this setting, Theorems \ref{thm:maindistr} and \ref{thm:mainfunction} 
yield the following result.

\begin{Thm} \label{cor:mainhyp}
Let the notation be as above. Let
$\Theta\subseteq \Pi$. For a hypergeometric
differential operator  $D_0$ let $D=\delta(x)^{2k} \cdot D_0$ be a
regularization of $D_0$. Then for every $u \in \dcsupptheta$ we have
\begin{equation*}
  \conv(\supp Du)=\conv(\supp u).
\end{equation*}
Equivalently, for every $f \in \ccisupptheta$, we have
\begin{equation*}
  \conv(\supp Df)=\conv(\supp f).
\end{equation*}
\end{Thm}

In the case of $\pedtheta{W}$-invariant functions, Corollary 
\ref{cor:mainhyp}
can be stated directly for the hypergeometric operators.

\begin{Thm} \label{cor:mainhypW}
Let the notation be as above. 
Let $D_0$ be a hypergeometric differential operator
and $f\in \cciwtheta$.
Then $D_0f \in \cciwtheta$ for all $f \in \cciwtheta$ and
\begin{equation*}
  \conv(\supp D_0 f)=\conv(\supp f)\, .
\end{equation*}
\end{Thm}
\pf As
$$\frak a=\bigcup_{w\in \widetilde{W}}\overline{w ({\frak a}^+)}=
\bigcup_{w\in \widetilde{W}/\pedtheta{W}}\overline{w({\frak a}_{\Theta})}$$
and $\supp (f)\subset {\frak a}_{\Theta}$ is compact
and $\pedtheta{W}$-invariant, it
follows that there exists a unique $\widetilde{f} \in
\cciwW$ so that $\widetilde{f}|_{\pedtheta{\frak a}}=f$. Let $D=
\delta^{2k} \cdot D_0$ be a regularization of $D_0$. Then 
\begin{equation}\label{eq:restrictions}
(D\widetilde{f})|_{\pedtheta{\frak a}}=D(\widetilde{f}|_{\pedtheta{\frak a}})
=Df\,.
\end{equation}
Suppose $D_0$ is the differential part of the Cherednik operator $T_0$
(see e.g. \cite{HeckBou}). Then $D_0\widetilde{f}=T_0\widetilde{f}$
by $\widetilde{W}$-invariance. Cherednik operators
map $\widetilde{W}$-invariant smooth functions into $\widetilde{W}$-invariant 
smooth functions. It follows
that $D_0\widetilde{f}$ is smooth and $\widetilde{W}$-invariant. 
>From  (\ref{eq:restrictions}) we therefore deduce that 
\begin{equation*}
  D_0 f=\frac{ (\delta^{2k} \cdot D_0\widetilde{f})|_{\pedtheta{\frak a}}}
{\delta^{2k}}
\end{equation*}
extends to be smooth and $\widetilde{W}$-invariant on 
$\pedtheta{\frak a}$. If $g$ is a continuous function, 
then $\delta^{2k}\cdot g$ and $g$ have the same support.
The Theorem therefore follows.
\qed

\subsection{Bessel differential operators}
\label{ex:bessel}
Let $G$ be a connected noncompact semisimple Lie group
with fintie center and $K\subset G$ a maximal
compact subgroup. Then $K=G^\theta$ for
some Cartan involution $\theta$. The Lie
algebra $\frak g$ decomposes into eigenspaces
of the derived homomorphism $\theta :\frak g \to
\frak g$:
$$\frak g=\frak k\oplus \frak p$$
where $\frak k$ is the $+1$ eigenspace
and $\frak p$ is the $-1$ eigenspace.
Consider $\frak p$ as a abelian Lie group. The
group $K$ acts linearly on $\frak p$ by $k\cdot X=\mathrm{Ad}(k)X$.
We can therefore consider the semi-direct product
$$G_0=\frak p \times_{\mathrm{Ad}} K\, .$$

Let $(\frak a, \Sigma, m)$ be as in Example \ref{ex:hyper}. 
The {\em Bessel differential operators} associated with $(\frak a, \Sigma, m)$
are the ``rational'' analogs of the hypergeometric differential operators.
They are $W$-invariant and the singularities of their meromorphic 
coefficients are cancelled by
multiplication be a power of the polynomial 
\begin{equation} \label{eq:pi}
\pi(x):=\prod_{\a \in \Sigma^+} \a(x)\,.
\end{equation}
When the triple $(\frak a, \Sigma, m)$ arises from the structure of
Riemannian symmetric space $G/K$ of the noncompact type, 
the Bessel differential operators
coincide with the radial parts on $\frak a^+$ 
of the $G_0$-invariant differential operators on the corresponding 
Riemannian symmetric space $G_0/K$ of Euclidean type.
For instance, the Bessel differential operator
\begin{equation} \label{eq:LO}
L_0:=L_{\frak a}+\sum_{\a \in \Sigma^+}m_\a \frac{1}{\a} \,\partial(y_\a)
\end{equation}
coincides with the radial part of the Laplace-Beltrami operator
on $G_0/K$.
The operator 
$L_0$ is the multivariable analog of the Bessel differential operator 
of Example \ref{ex:diffrankone}.
The Bessel differential operators can be constructed
as ``differential parts'' of differential-reflection operators associated
with $W$-invariant polynomials functions $p$ on $\frakacs$ by means
of the Dunkl operators \cite{deJeuPW}.
Suppose $D_0$ is the Bessel differential operator associated with $p$. 
Let $k \in \N$ be chosen so that $D:=\pi(x)^{2k} D_0$ has real analytic 
coefficients. Then the principal symbol of $D$ is
\begin{equation*}
  \sigma(D)(x,\l)=p_h(\l) \pi(x)^{2k}=p_h(\l)\pedtheta{\pi}(x)^{2k}
P(\l,x)\,,
\end{equation*}
where $p_h(\l)$ is the highest homogeneous part of $p(\l)$
and $P(\l,x)=\prod_{\a \in \Sigma^+ \setminus \rootstheta^+} \a(x)^{2k}$
is a homogeneous $\pedtheta{\frak a}$-elliptic polynomial of degree $0$
in $\l$. As for the hypergeometric differential operators, we obtain
the following corollary of Theorems \ref{thm:maindistr} and 
\ref{thm:mainfunction}.

\begin{Thm} \label{cor:mainbessel}
Let the notation be
as above. Let $\Theta$ be a fixed set of
positive simple roots in $\Sigma$. For
a Bessel differential operator $D_0$ let $D=\pi(x)^{2k} \cdot D_0$ be a
regularization of $D_0$. Then, for $u \in \dcsupptheta$ we have
\begin{equation*}
  \conv(\supp Du)=\conv(\supp u).
\end{equation*}
Equivalently, for every $f \in \ccisupptheta$, we have
\begin{equation*}
  \conv(\supp Df)=\conv(\supp f).
\end{equation*}
If $f \in \cciwtheta$, then $D_0f \in \cciwtheta$
and
\begin{equation*}
  \conv(\supp D_0 f)=\conv(\supp f)\, .
\end{equation*}
\end{Thm}

\subsection{Shift operators}\label{ex:shift}
For the study of hypergeometric differential operators corresponding to
different multiplicity functions (see Example \ref{ex:hyper}), 
Opdam introduced certain $W$-invariant differential operators,
called \emph{shift operators}.  For simplicity,
we treat here only the case of a reduced root system $\Sigma$ and 
we refer the reader to \cite{HS}, Part I, Ch. 3 for the general case.
As in the case of hypergeometric differential operators, Opdam's
shift operators have meromorphic coefficients with singularities along 
the root hyperplanes, and their singularities are cancelled by 
multiplication by powers of $\Delta$.
A \emph{shift} is a multiplicity function $l$ which is even 
(i.e. such that $l_\a \in 2\Z$ for all $\a \in \Sigma$). 
Opdam's shift operators exists for any shift $l$. They can be considered 
as generalizations of the hypergeometric differential operators, which form 
the space $\mathcal S(0)$ of shift operators of shift $0$. 
For an arbitrary nonzero shift $l$, the space 
$\mathcal S(l)$ of shift operators of shift $l$ is a free rank-one (right)
$\mathcal S(0)$-module generated by a certain shift operator $G(l)$.
See \emph{loc. cit.}, Theorem 3.3.7. Suppose furthermore that either
$l_\a>0$ for all $\a$ or $l_\a<0$ for all $\a$. Then $G(l)$ is obtained 
by composition of \emph{fundamental shift operators} $G_{S,\pm}$, shifting 
of $\pm 2$ the multiplicities along each Weyl group orbit (\emph{loc. cit.}, 
p. 42). The principal symbol of each fundamental shift operator is 
\begin{equation*}
  \sigma(G_{S,\pm})(x,\l)=c \cdot \Delta_S(x)^{\pm 1} \prod_{\a \in S^+} 
  \l(x_\a)\,,
\end{equation*}
where $S^+=\Sigma^+ \cap S$ and $\Delta_S(x)=\prod_{\a \in S^+} \sinh \a(x)$. 
See \emph{loc. cit.}, Remark 3.3.8. The analysis done in Example
\ref{ex:hyper} for the hypergeometric differential operators, in particular 
Corollary \ref{cor:mainhyp}, easily extends also
to the fundamental shift operators, and hence to all shift operators.

A similar argument can be also applied to the shift operators
associated with the Bessel differential operators.
For more information on the latter we refer to 
\cite{HeckBirk}, where they are studied in the general case in which $W$ is an 
arbitrary finite Coxeter group. 

\subsection{Hamiltonian systems}
A wide class of integrable Hamiltonian systems associated with root systems 
were introduced by Olshanetsky and Perelemov in \cite{OPe76}.
Let $(\frak a,\Sigma,g)$ be a triple consisting of a finite dimensional
Euclidean space $\frak a$, a root system $\Sigma$ in $\frak a^*$, and a
real-valued multiplicity function $g$ on $\Sigma$.   
These integrable systems are described by a Hamiltonian of the form
\begin{equation*}
  H=-\frac{1}{2} \inner{p}{p}+U(q), \qquad p,q \in \frak a\,,
\end{equation*}
with potential energy 
\begin{equation*}
  U(q):=\sum_{\a \in \Sigma^+} g_\a^2 v(\inner{q}{\a})
\end{equation*}
where the function $v$ has five possible forms. Here we consider the cases
I, II and V as listed in \cite{Pe90}, (3.1.14) and (3.8.3):
\begin{alignat*}{2}
&\mathrm{(I)} \qquad &v(\xi)&=\xi^{-2}\,,\\
&\mathrm{(II)} \qquad &v(\xi)&=\sinh^{-2}\xi\,,\\
&\mathrm{(V)} \qquad &v(\xi)&=\xi^{-2}+\omega^2 \xi^2\,.
\end{alignat*}
The Calogero Hamiltonian of Example \ref{ex:Calogero} is a special instance
of case I. The associated Schr\"odinger operators are the linear 
partial differential operators
\begin{equation*}
  S=-\frac{1}{2} L_{\frak a}+U(x),\qquad  x \in \frak a. 
\end{equation*}
They are invariant with respect to the Weyl group $\widetilde{W}$ of
$\Sigma$. The operator $S$ can be regularized by multiplication by the 
polynomial $\pi(x)$ of (\ref{eq:pi}) in the cases I and V, and by 
multiplication by $\delta$
as in (\ref{eq:delta}).
The possibility of applying Theorems \ref{thm:maindistr} and 
\ref{thm:mainfunction} to these regularized differential operators can 
be proven as in Examples \ref{ex:hyper} and \ref{ex:bessel}.
In fact, there is a close relation between
the operator $S$ of case II, 
resp. of case I, and the operator $L$ of 
(\ref{eq:L}), resp. to the operator $L_0$ of (\ref{eq:LO}). See e.g. 
 \cite{HS}, Theorem 2.1.1, and \cite{OPe76}. 

For an overview on the role
of the Hamiltonian systems treated in this example in different
areas of theoretical physics and mathematics, we refer the reader 
to \cite{Dj00}.


\section{Applications to symmetric spaces}
\label{section:appl}
 
\noindent
In this section we apply Theorem
\ref{thm:mainfunction} to differential
operators on symmetric spaces. In
particular we give a new proof of the
$D$-convexity of the Riemannian symmetric
space $G/K$, c.f. \cite{He73,vdBSconv}. Our
proof uses Theorem \ref{thm:mainfunction},
applied to the radial part of invariant
differential operators and is,
as far as we know, new. Our notation is the same
as in Subsection \ref{ex:bessel}.

Let us start by recalling the notation
from Subsection \ref{ex:bessel}. Here
$G$ is a connected noncompact semisimple
Lie group with 
finite center. Furthermore $\theta :G\to G$ is
a Cartan involution and $K=G^\theta$ the
corresponding maximal compact subgroup.
Denote by $\kappa : G\to G/K$, the
natural projection $g\mapsto gK$. We
denote also by $\theta$ the derived involution
on $\frak g$. Then $\frak g =\frak k\oplus \frak p$,
where $\frak k$ is the $+1$-eigenspace of $\theta$, and
$\frak p$ the $-1$-eigenspace. Then $\frak k$ is
the Lie algebra of $K$.
Fix a Cartan subspace $\mathfrak{a}$, that is a
maximal abelian subspace  of $\mathfrak{p}$.
The Killing form on $\frak g$ defines an Euclidean inner product on 
$\frak a$. For $\a \in \frak a^*$ set 
$\mathfrak{g}^{\alpha}=\{y\in \mathfrak{g}^{\alpha}: 
\text{$[x,y]=\a(x) y$ for all $x \in \frak a$}\,\}$.
Then the set $\Sigma$ consisting of all $\a \in \frak a^* \setminus \{0\}$
for which $\mathfrak{g}^{\alpha}\neq \{0\}$ is a (generally non-reduced) 
 root system as defined in Section \ref{section:examples}. 
It is called the \emph{(restricted) root system} 
of $(\mathfrak{g}, \mathfrak{a})$. The multiplicity $m_\a$ of $\a\in \Sigma$
is defined as the dimension of $\mathfrak{g}^{\alpha}$. The map 
$m$ given by $m(\a):=m_\a$ is a multiplicity function on $\Sigma$.
This construction associates with the Riemannian symmetric space $G/K$ 
a triple $(\mathfrak{a}, \Sigma, m)$.   As before
we denote by $\widetilde{W}$ the corresponding
Weyl group.
 
Let $\Sigma^+$ be a choice of positive roots and let 
$\mathfrak{a}^+:=\{x \in \frak a: \text{$\a(x)>0$ for all 
$\a \in \Sigma^+$}\}$ be the corresponding positive chamber.
Denote by $\exp:\mathfrak{g} \to G$ the exponential map. Then
$A:=\exp\frak a$ is an abelian subgroup of $G$ diffeomorphic to $\frak a$.
We set $A^+:=\exp \frak a^+$.  The map $K\times A\times K\ni (k_1,a,k_2)
\mapsto k_1ak_2\in G$ is surjective and the $A$-component
is unique up to conjugation by an element
of $\widetilde{W}$.
Hence every $K$-bi-invariant subset of $G$ is of the form
$K(\exp B) K$ where $B$ is a $\widetilde{W}$-invariant subset of $\frak a$.
Moreover, $B$ is compact if and only if $\kappa (B)\subset G/K$ is.
Let $C^\infty_c(G/K)$ denote the space of compactly supported
smooth functions on $G/K$. Using the map $C_c^\infty (G/K)\to C_c^\infty (G/K)^K$,
$f\mapsto f\circ \kappa$, we will often identify 
smooth functions on $G/K$ with the
corresponding right $K$-invariant function on $G$.

The decomposition $G=KAK$ yields the following lemma.

\begin{Lemma} \label{lemma:suppfKAK}
Let $f \in C_c^\infty(G/K)$ and suppose that $\supp f$ is $K$-bi-invariant.
Then $\supp f=K (\supp f|_A)K$ where
$f|_A$ denotes the restriction of $f$ to $A$. Moreover, $\supp f|_A$
is a $\widetilde{W}$-invariant compact subset of $A$. 
\end{Lemma}

Denote by $\D (G/K)$ the (commutative) algebra of $G$-invariant
differential operators on $G/K$.
We identify  $A^+$ can be identified with the
submanifold $\kappa (A^+)\subset G/K$. Then, for every $D \in \D(G/K)$,
there is a unique $\widetilde{W}$-invariant differential operator,
$\omega(D)$ on $A^+$,
called the \emph{radial part of $D$}, so that for all $f \in C^\infty(G/K)$ one has
\begin{equation} \label{eq:radialpart}
  (Df)|_{A^+}=\omega(D) (f|_{A^+})
\end{equation}
See \cite{He2}, p. 259.

Define a differential operator on $\frak a$, also
denoted by $\omega (D)$, by
\begin{equation*}
  \omega(D)g := \omega(D)(g \circ \exp^{-1}) \circ \exp\,,
\qquad  g \in C^\infty(\frak a)\,.
\end{equation*}
In this way, we can consider $\omega(D)$ as a (singular) $\widetilde{W}$-invariant
differential operator on the Euclidean space $\frak a$.
To simplify our notation, we shall adopt the identification $A\equiv 
\frak a$ using the exponential map. Then
$\exp:=\id$. We then write $f|_{\frak a}$
instead of $f|_A$, and the above mentioned decomposition of a 
$K$-bi-invariant subset of $G$ will be written 
as $KBK$ instead of $K(\exp B)K$.
With these identifications, the operator $\omega(D)$ turns out to be a 
hypergeometric differential operator associated with the triple
$(\frak a,\Sigma,m)$ as considered in Subsection \ref{ex:hyper}. When
$f \in C^\infty_c(G/K)$ is $K$-invariant, we can moreover extend
(\ref{eq:radialpart}) by $\widetilde{W}$-invariance to obtain
\begin{equation} \label{eq:radialpartona}
  (Df)|_{\frak a}=\omega(D) (f|_{\frak a})
\end{equation}

\begin{Lemma}\label{lemma:almostDconvex}
Let $G/K$ be a Riemannian symmetric space of the noncompact type. 
 Let $D \in \D(G/K)$ and let $B$ be a compact, convex and 
$\widetilde{W}$-invariant subset of $\frak a$. Set $X_B:=KBK$.
Then for all $f \in C_c^\infty(G/K)$, we have
$$\supp f \subseteq X_B \quad\text{if and only if}\quad \supp(Df) \subseteq X_B.$$
\end{Lemma}
\pf
We need to prove that $\supp(Df) \subseteq X_B$ implies $\supp f \subseteq X_B$.
The first step is, as in \cite{vdBSconv}, to reduce  to the case
in which the support of $f$ is left-$K$-invariant. The function $f$ 
can in fact be expanded as a sum of $K$-finite functions. 
Since $X_B$ is left-$K$-invariant, we will obtain $\supp f \subseteq X_B$ 
if the support of each $K$-finite summand is contained in $X_B$. 
We can therefore assume that $f$, and hence $Df$, are
$K$-finite. As the support of a  $K$-finite function is left $K$-invariant,
we obtain from Lemma \ref{lemma:suppfKAK} that $\supp f=
K(\supp f|_{\frak a})K$ and $\supp(Df)=
K\big(\supp(Df)|_{\frak a}\big)K$ where $\supp f|_{\frak a}$ and 
$\supp(Df)|_{\frak a}$ are $\widetilde{W}$-invariant and compact. 
Since $\supp(Df) \subseteq X_B$, we then conclude that
$\supp(Df)|_{\frak a} \subseteq B$. Since $\omega(D)$ is a hypergeometric
differential operator, Corollary \ref{cor:mainhypW} with $\Theta=\Pi$
yields that $\conv(\supp f|_{\frak a})=\conv(\omega(D)\supp f|_{\frak a})=
\conv(\supp(Df)|_{\frak a}) \subseteq B$.
This proves the required inclusion.
\qed

\begin{Thm}\label{thm:Dconvex}
Let $G/K$ and $D \in \D(G/K)$ be as in Lemma \ref{lemma:almostDconvex}. Then 
$G/K$ is $D$-convex, that is for every compact subset $S$ of 
$G/K$ there is a compact
set $S'$ so that for every $f \in \mathcal C^\infty_c(G/K)$ the 
inclusion $\supp(D^t f) \subseteq S$ implies  $\supp f \subseteq S'$. 
Here $D^t$ denotes the formal transpose of $D$. 
\end{Thm}
\pf
Choose $S_1 \subseteq G$ so that $S_1/K=S$
The set $KS_1H$ is $K$-bi-invariant, hence $KS_1K=X_B:=KBK$
for some $\widetilde{W}$-invariant compact subset $B$ of $A \equiv \frak a$.
Since $D^t \in \D(G/K)$ we can apply Lemma \ref{lemma:suppfKAK} to it.
So $\supp(f) \subseteq X_{\conv B}$. We can then select $S':=X_{\conv B}/K$.
\qed

We conclude this section by a short
discussion of the application of the theorem
of our support theorem in harmonic analysis corresponding
the the $\Theta$-hypergeometric
transform, c.f. \cite{OP4}. Because of the technical nature of this
application, a detailed exposition would require a certain amount of 
notation and of background information. Instead, we prefer here just to 
outline the main ideas involved. We refer the interested reader to 
\cite{OP4} for further information.  

As already remarked in the introduction, the support theorem plays
a role in the study of the \textit{Paley-Wiener space}.
This is the set of images under a suitable generalization $\mathcal F$ 
of the Fourier transform, of the compactly supported smooth functions.
The Paley-Wiener space generally consists of 
entire or meromorphic functions with exponential growth and possibly 
satisfying additional symmetry conditions.
The size the support of the original function is linked to the exponential 
growth of its Fourier transform. The thrust of Paley-Wiener type theorems 
is usually to prove that, if a function $g$ in the Paley-Wiener space
has a given exponential growth, then the support of the associated 
``wave packet'' 
$\mathcal I g$ (obtained by formal application to $g$ of the inverse 
Fourier transform $\mathcal I$ of $\mathcal F$) 
is compact and has the correct size.
In the classical situation of the Fourier transform on $\R^n$, this is proven 
by a suitable ``shift'' of contour of integration by means of Cauchy's theorem.
The shift is allowed because $\mathcal I g$ is given by integration of
an entire function of exponential type and rapidly decreasing. 
A suitable generalization
of this argument was applied also to the spherical Fourier transform on 
Riemannian symmetric spaces. See \cite{He2}, Ch. IV, \S 7.2, or \cite{GV}, 
\S 6.6. 
For pseudo-Riemannian symmetric spaces the situation is more complicated.
Here the wave packages are integrals of meromorphic functions and the required
shift of integration would generally require addition of certain 
``residues''. For some kinds of Fourier transforms, like the spherical
Fourier-Laplace transform on noncompactly causal symmetric spaces
with even multiplicities \cite{OP4}, the singularities of the integrand are
cancelled by multiplication by
a certain polynomial function $p$. Hence, for every function $g$ in the
Paley-Wiener space the wave-packet $\mathcal I(pg)$ is given by integration
of an entire function. Furthermore, there is an invariant differential
operator $D$ so that $D(\mathcal Ig)=\mathcal I(pg)$. The possibility of
shifting the contour of integration in $\mathcal I(pg)$ allows us then
to determine the size of the support of $D(\mathcal Ig)$. The generalization
of the support theorem of Lions and Titchmarsch presented in this papers
allows us finally to recover the size of the support of $\mathcal Ig$.
A concrete application
of this procedure using Theorem \ref{thm:mainfunction} can be found 
in \cite{OP4} in the context of the $\Theta$-hypergeometric transform.
The latter transform, stated in a setting of transforms associated with 
root systems, is a common generalization of Opdam's hypergeometrical 
transform (hence of Harish-Chandra's spherical transform on Riemannian 
symmetric spaces of noncompact type) and of the Fourier-Laplace transform 
on noncompactly causal symmetric spaces \cite{HO}.
A similar application would also be needed for determining the Paley-Wiener
space for the transform associated with the $\Theta$-Bessel functions
of \cite{BSO04}.

\end{document}